\documentclass[12pt]{amsart}

\usepackage{amscd}
\usepackage{amsmath,amsfonts,amsthm,epsfig,latexsym,graphicx,amssymb,psfrag}
\usepackage[english]{babel}

\newtheorem{prop}{Proposition}[section]

\newtheorem*{theor}{Theorem}

\newcommand{\pf}{{\bf Proof}: }

\def\R{I\kern -0.37 em R}
\def\N{I\kern -0.37 em N}
\def\Z{I\kern -0.37 em Z}

\def\supess_#1{\mathop{\rm supess}\limits_{#1}}
\def\infess_#1{\mathop{\rm infess}\limits_{#1}}

 \def\RR{{\mathbb R}}

\begin{document}

\title[$C^1$-BS action on $S^1$.]{\bf $C^1$- Actions of   Baumslag-Solitar groups   on
  $S^1$.}
\author{Nancy Guelman and Isabelle Liousse} \thanks{This paper was partially
supported by Universit\'{e} de Lille 1, PEDECIBA, Universidad de la
Rep\'{u}blica and the PREMER project.}

\address{{\bf  Nancy Guelman}
IMERL, Facultad de Ingenier\'{\i}a, Universidad de la Rep\'ublica,
C.C. 30, Montevideo, Uruguay.  \emph{nguelman@fing.edu.uy}, }

\address{{\bf    Isabelle Liousse}, UMR CNRS 8524, Universit\'{e} de Lille1,
59655 Villeneuve d'Ascq C\'{e}dex,   France.  \emph {liousse@math.univ-lille1.fr},  }

\begin{abstract} Let $BS(1, n) =< a, b \ | \ aba^{-1} = b^n >$  be the solvable
  Baumslag-Solitar group, where $ n\geq 2$. It is known that $BS(1, n)$
  is isomorphic to the group generated by the two affine maps  of the line :
  $f_0(x) = x + 1$ and $h_0(x) = nx $. The action  on $S^1 = \RR \cup {\infty}$    generated by these two affine maps $f_0$ and $h_0 $ is called the standard affine one.
  We prove that any representation of $BS(1,n)$ into $Diff^1(S^1)$ is (up to a finite index subgroup) semiconjugated to the standard affine action.

\end{abstract}

\maketitle
\date{}

\section {Introduction.}

\medskip

This paper is about the dynamics of the actions of the solvable
Baumslag-Solitar group --$BS(1, n) =< a, b \ | \ aba^{-1} = b^n >$, where $ n\geq
2$-- on the circle $S^1$.

\medskip

It is  known that  $BS(1,n)$  has many  actions on $S^1$.
 The standard  action on $S^1 = \RR \cup {\infty}$ is
 the action   generated by the two affine maps $f_0(x) = x + 1$ and $h_0(x) = nx $ (where $f_0 \equiv b$ and $h_0 \equiv a $).

Many people has studied actions of solvable groups on one-manifolds, for example  Plante \cite{P}, Ghys \cite{Gh}, Navas \cite{Na}, Farb and Franks\cite{FF}, Moriyama\cite{Mo} and Revelo and Silva\cite{RS}.

There exist some results concerning to  $BS(1,n)$-action on $S^1$. Many of them are scattered in different articles and some of them are not so easy to find, so our aim is to present the state of the art for the case of  the action of $BS(1,n)$-group on $S^1$, in the case $C^1$.

It is known that any $C^2$ $BS(1,n)$-action on $S^1$ admits a finite
orbit. This fact was proven by Burslem-Wilkinson in  \cite{BW}. In
fact they gave a classification (up to
conjugacy) of representations  $ \rho: BS(1,n)
\rightarrow Diff^r(S^1)$ with $r\geq 2$ or $r=\omega$.

%

These results are proved by using a dynamical approach. The dynamics of $C^2$
$BS(1,n)$-actions on $S^1$ is now well understood, due to Navas work on
solvable group of circle diffeomorphisms (see \cite{Na}).
We will prove the same statement  of  \cite{BW} but in the case that the action is $C^1$, that is,  any  $C^1$ $BS(1,n)$-action on $S^1$ admits a finite orbit.

\bigskip

Our aim is giving a classification of faithful actions  of $BS(1,n)$ on
   $S^1$, in the case $C^1$.

So, our main result is the following:

\begin{theor}
Any $C^1- BS(1,n)$ action on $S^1$ preserving orientation, $<f, h>$, is (up to a finite index subgroup $<f^{n-1}, h^{m}>$) semiconjugated to the standard
affine action.
\end{theor}
The proof uses in a crucial way (a slightly extended version of) an argument due to Cantwell and Conlon (\cite{CC}), and its use was suggested to us by  A. Navas. It should be pointed out that, using this argument, Navas has recently obtained a counter-example to the converse of Thurston Stability Theorem (see \cite{Nava}), and it is very likely that Cantwell-Conlon argument will still reveal very fruitful for obtaining new obstructions to $C^1$ actions on 1-dimensional manifolds

   \section {Existence of a global finite orbit.}

 As before, let
$<f,h>$ be a faithful action  of $BS(1,n)$ on
   $S^1$, where  $h\circ f\circ h ^{-1} = f^n$.

   The aim of this section is proving the existence of a common periodic point of $f$ and $h$.

We begin by proving that the rotation number of $f$ is rational.

\begin{prop}Let $f,h: S^1 \rightarrow S^1$ be  homeomorphisms preserving
orientation  verifying $h\circ f\circ h ^{-1} = f^n$. Then there exists $ l \in \mathbb{N}$ such that the
rotation number  of $f$, $\rho (f)=\frac{l}{n-1}$.

\end{prop}

\pf

Since $h\circ f\circ h ^{-1} = f^n$, we have that there exists $l \in \mathbb{Z}$ such that $$n \rho (f)= \rho(f) + l,$$ then 
$\rho (f)=l/(n-1)$. $\blacksquare$

We will prove that not only the rotation number of $f$ is rational
but also the rotation number of $h$ it is.
\begin{prop}Let $f,h: S^1 \rightarrow S^1$ be $C^1$- diffeomorphisms preserving
orientation  verifying $h\circ f\circ h ^{-1} = f^n$. If the
$BS(1,n)-$ action on $S^1$, $<f,h>$, is faithful then $\rho (h)$,
the rotation number of $h$,  is rational. \label{phoh}
\end{prop}
We suppose that $\rho (h)$ is irrational, then we have two cases:
\begin{enumerate}
\item $h$ is conjugated to an irrational rotation.
Note that the periodic points of $f$ are preserved by $h$: let $q$
such that $f^k(q)=q$, then $ f^{nk}(h(q))=h(f^k(q))=h(q)$. So $h(q)$
is a  periodic point for $f$.

 It follows that the periodic points of
$f$ are dense in $S^1$. This implies that there exists $m$ such that
$ f^m = Id$ contradicting that the action is faithful.

\item The minimal set of $h$ is a Cantor set, $K$. Then there exists
an arc $J\subset K^c \subset S^1$ fixed by $f$ (its end points are
fixed) where $f|_J \neq Id$ and length of $J$ is very small. Let
$\epsilon >0$ verifying $(1-\epsilon)^2 >\frac{3}{4}$, we choose $J$
sufficiently small  in order to $length (h^{-s}(J))$  (it is possible since $J$ is an $h$-wandering interval) is small enough   to guarantee $f'\geq 1-\epsilon$ for any $x\in J$ and
$\frac{h'(x)}{h'(y)}\geq 1-\epsilon$ for any $x, y \in J$ and also  $f'\geq 1-\epsilon$ for any $x\in h^{-s}(J)$ and
$\frac{h'(x)}{h'(y)}\geq 1-\epsilon$ for any $x, y \in h^{-s}(J)$, for any $s>0$.

 Let
$x\in J$ and $ I\subset J$ be the open arc between $x$ and $f(x)$. It is
easy to see that $\cup_{k \in \mathbb{Z}}f^k(I) \subset J$ and
$f^k(I) \cap f^j(I)= \emptyset$ if $k \neq j$. Let us define the map
$\Psi(\alpha_0,\alpha_1,\alpha_2,..., \alpha_m) $ with $\alpha_i \in
\{0,1\}$
 as the map $$ (h^m \circ f^{\alpha_m} \circ h^{-m}) \circ ...
 \circ (h^2 \circ f^{\alpha_2} \circ h^{-2}) \circ (h \circ f^{\alpha_1} \circ
 h^{-1})\circ f^{\alpha_0}.$$
 Since $h^i \circ f^{\alpha_i} \circ h^{-i} = f^{n^i}$ if
 $\alpha_i=1$ and $h^i \circ f^{\alpha_i} \circ h^{-i} =Id$,
 otherwise; it follows that
 $\Psi(\alpha_0,\alpha_1,\alpha_2,..., \alpha_m)(I)= f^{\beta} (I) $
 where $ \beta = \sum_{i | \alpha_i=1} n^i$.

 Therefore, for $(\alpha_0,\alpha_1,\alpha_2,...,
 \alpha_m) \neq (\alpha'_0,\alpha'_1,\alpha'_2,..., \alpha'_m)$ it
 holds that $\Psi(\alpha_0,\alpha_1,\alpha_2,..., \alpha_m)(I) \cap \Psi(\alpha'_0,\alpha'_1,\alpha'_2,...,
 \alpha'_m)(I)$ is empty.
 Hence, $$ \mathcal{I}=\bigcup_{\{0,1\}^{m+1}}\Psi(\alpha _0,\alpha _1,\alpha _2,..., \alpha _m)(I)$$
  is the union of $2^{m+1}$ disjoint arcs included in $J$.

We claim that the length of any arc $\Psi(\alpha _0,\alpha _1,\alpha
_2,..., \alpha _m)(I)$ may be lower bounded by $|I|
(\frac{3}{4})^{m+1}$  :

Notice that $$\Psi(\alpha_0,\alpha_1,\alpha_2,..., \alpha_m)=h^k
\circ f \circ h^{-l_1}\circ f \circ h^{-l_2}\circ f \circ...\circ f
\circ h^{-l_r} $$ where $l_1+l_2+...+l_r=k$, $l_i>0$ for $i=1,...,
r-1$ and $l_r\geq 0$. Hence, $$\Psi'(\alpha_0,\alpha_1,\alpha_2,...,
\alpha_m)(u)=\prod_{i=1}^r f'(x_i) \prod_{j=1}^k
\frac{h'(y_j)}{h'(w_j)},$$ where $x_i, y_j$ and $w_j$ are well
defined points, $x_i \in \cup_{s \in \mathbb{N}} h^{-s}(J)$ and $y_j, w_j \in h^{-j}(J)$.

Therefore, there exist points $\widehat{x_i}$  for $i=1,..,r$ , $\widehat{y_j}$ and $\widehat{w_j}$  for $j=1,..,k$ such that
the length of
$\Psi(\alpha_0,\alpha_1,\alpha_2,..., \alpha_m)(I)$,
$$|\Psi(\alpha_0,\alpha_1,\alpha_2,...,
\alpha_m)(I)|=|I|\prod_{i=1}^r |f'(\widehat{x_i}) |\prod_{j=1}^k
\frac{|h'(\widehat{y_j})|}{|h'(\widehat{w_j})|}\geq|I| (1-\epsilon)^{r+k}.$$ Since
$r\leq k\leq m+1$ it follows that
$$|\Psi(\alpha_0,\alpha_1,\alpha_2,..., \alpha_m)(I)| \geq |I|
(1-\epsilon)^{2m+2}\geq |I| (\frac{3}{4})^{m+1}.$$

Then, the length of $\mathcal{I} \geq 2^{m+1} |I|
(\frac{3}{4})^{m+1}$ which tends to $\infty$ when $m\rightarrow
\infty$. This is a contradiction with $\mathcal{I}\subset J$.

Therefore, the minimal set of $h$ is not a Cantor set.

\end{enumerate}
Hence, we have proved that $\rho(h)$ is rational.$\blacksquare$

\begin{prop}There exists $m\in \N$ such that $f^{n-1}$ and $h^m$ have a
common fixed point
\end{prop}
\pf Let $m$ such that $\rho (h^m)=0$. Let $q$ be a fixed point for
$f^{n-1}$ (we have already seen that there exists a fixed point for
$f^{n-1}$), then $\{ h^{lm}(q)\}$ is included in the set of fixed points of
$f^{n-1}$  since $h$  preserves the fixed points of $f^{n-1}$.

Let $u$ be an accumulation point of $\{ h^{lm}(q)\}$.
 It follows from continuity of $f$ that $u$ is a fixed point for $f^{n-1}$
 and since $\Omega (h^m)$ is included in the set of fixed points of $h^m$ then $u$ is also a fixed point for $h^m$.$\blacksquare$

\vspace{2mm}

Since the $C^1$ diffeomorphisms $f^{n-1}$ and $h^m$ have a common fixed
point, we can study $<f^{n-1}, h^m>$, a $C^1$-action of $BS(1, (n-1)n^m)$ on
the interval $[0,1]$ instead of $S^1$.

\section{Semiconjugation to the standard affine action.}

   Recall that the standard  affine action on $S^1 = \RR \cup {\infty}$ is
 the action   generated by the two affine maps $f_0(x) = x + 1$ and $h_0(x) = nx $.

\vspace{3mm} Following results or ideas of Cantwell-Conlon, Navas
and Rivas, from now on we will prove that any faithful $C^1- BS(1,n)$ action
on $S^1$ is semiconjugated to the standard affine action.

\vspace{2mm} Due to a classical result that appears, for example in
``Groups acting on the circle" by E. Ghys ( see \cite{Gh}), it is
known that for a countable infinite group $G$, $G$ is left orderable
if and only if $G$ acts faithfully on the real line by orientation
preserving homeomorphisms. Let us note that for proving that $G$
acts faithfully on the real line E. Ghys constructed the
\textbf{dynamical realization of a left ordering}. For this
construction he fixed an enumeration $\{g_i\}$ of $G$. He defined an
order preserving map $t: G \rightarrow \mathbb{R}$, in such a way
that $ g(t(g_i))=t(g g_i)$. This action was extended to the closure
of $t(G)$, and later to the whole real line. As Rivas noted (see
remark 4.4 of \cite{Ri}), and it was proved by Navas (see
\cite{Nav}) that the dynamical realization associated to differents
enumerations of $G$ (but the same ordering) are topological
conjugated.

 It was proven by Rivas (see \cite{Ri})) that the set of left orderings of $BS(1,n)$ is made up of four Conradian orderings and an uncountable set of non Conradian left orderings. Each one of these infinetely many non Conradian orderings can be realized as an induced ordering that comes from an affine action on $BS(1,n)$. Moreover, in the proof of this result  it was shown that the dynamical realization of any non Conradian ordering is semiconjugated to the standard affine action. (In fact, Rivas proved this result for $BS(1,2)$ but the proof for $BS(1,n)$ is the same).

 As Ghys and Rivas results are ``topological", they hold in an interval instead of the real line.

We will call an ``exotic" action to the $BS(1,n) $ one that is induced by  a Conradian ordering.

Let $\ll f\gg $ be the largest abelian subgroup containing $f$. For any  ``exotic" action in the interval $[0,1]$ (that is induced by one of the four  Conradian orderings),  $\ll f\gg $  is a convex subgroup (in the sense of ordening, see the proof of Proposition $4.1$ of  \cite{Ri}).
 Then, for a point $x_0$ in $[0,1]$, the sequence $\{ f^n(x_0)\}$ is lower and upper bounded, then there exist
 $\lim_{n \rightarrow \infty}f^n(x_0)=d<1$ and $\lim_{n \rightarrow -\infty}f^n(x_0)=c>0$ two $f$- fixed points verifying that there is no  fixed point for $f$ in $(c,d)$. Moreover, $h(c,d)$ is disjoint of $(c,d)$.


 %
%

It was proven by   Cantwell and Conlon
(\cite{CC}) that any ``exotic"  $BS$- action on the interval $I$ can not be
$C^1$. In fact, the proof of  item (2) of Proposition \ref{phoh}
follows Cantwell-Conlon's proof.

It follows that any $C^1$ $BS$- action on an interval $I$ is
semiconjugated to the standard affine action. The same holds for
$C^1$ $BS$- action  on $S^1$.

So, we have proved

\begin{theor}
Any $C^1- BS(1,n)$ action on $S^1$ is (up to a finite index subgroup $<f^{n-1}, h^{m}>$) semiconjugated to the standard
affine action.
\end{theor}


\begin{thebibliography}{99}





\bibitem{BW} Burslem L.,Wilkinson A. \emph{Global rigidity of solvable group action on $S^1$.} Geom. Topol. 8 (2004) 877-924.


\bibitem{CC} Cantwell J., Conlon L. {\em An interesting class if $C^1$ foliations.} Topol. and its Applic. 126 (2002) 281-297.

\bibitem{FF} Farb B. and Franks, J.  \emph{Groups of homeomorphisms of one manifolds I: Actions of nonlinear groups.} Preprint, (2001).

\bibitem{Gh} Ghys E., \emph{Groups acting on the circle.} Enseign. Math. (2) 47, no 3-4, (2001),329-407.

\bibitem{Mo} Moriyama Y. \emph{Polycyclic groups of diffeomorphisms on the half-line.} Hokkaido Math. Jour.,23, no 3 (1994) 399-422.
\bibitem{Na} Navas A.{ \em Groupes r\'{e}solubles de diff\'{e}omorphismes de l' intervalle, du circle et de la droit.} Bull. Bras. Math. Soc. N.S. 35, no 1, (2004) 13-50.

\bibitem{Nav} Navas A.{ \em On the dynamics of left orderable groups.}  Ann. Inst. Fourier
(to appear)
\bibitem{Nava} Navas A.{ \em A finitely generated, locally indicable group with no faithful action by $C^1$ diffeomorphisms of the interval}  Geom. Topol. 14 (2010) 573– 584

\bibitem{P} Plante,J. {\em  .}  (2009), .
\bibitem{RS} Revelo J. and Silva R. \emph{The multiple ergodicity of nondiscrete subgroups of $Diff^w(S^1)$.} Mosc. Math. J. 3, (2003),123-171.
\bibitem{Ri} Rivas C. {\em  On spaces of Conradian group ordenings.} J. Group Theory (2009), 1-17.
\end{thebibliography}
\end{document}